\documentclass[11pt]{pjm}
\usepackage{amssymb,latexsym}
\usepackage{graphicx}
\usepackage{epstopdf}

\renewcommand{\int}{\operatorname{int}\!}

\newcommand{\Arf}{\operatorname{Arf}}

\newcommand{\Z}{\mathbb{Z}}
\newcommand{\N}{\mathbb{N}}

\newcommand{\cW}{\mathcal{W}}

\def\theenumi{\roman{enumi}}
\def\labelenumi{(\theenumi)}

\title[Whitney towers, gropes and the Arf invariant]{Simple Whitney towers,
half-gropes and the Arf invariant of a knot}

\author[Schneiderman]{Rob Schneiderman}
\address{Courant Institute of Mathematical Sciences,
New York University\\
251 Mercer Street New York NY 10012-1185 USA}
\email{schneiderman@courant.nyu.edu}
\newtheorem{thm}{Theorem}
\newtheorem{lem}[thm]{Lemma}

\theoremstyle{definition}

\newtheorem{defi}{Definition}

\theoremstyle{remark}

\begin{document}

\begin{abstract}

A geometric characterization of the Arf invariant of a knot in the
3--sphere is given in terms of two kinds of 4--dimensional
bordisms, {\em half-gropes} and {\em Whitney towers}. These types
of bordisms have associated complexities {\em class} and {\em
order} which filter the condition of bordism by an embedded
annulus, i.e. knot concordance, and it is shown constructively
that the Arf invariant is exactly the obstruction to cobording
pairs of knots by half-gropes and Whitney towers of arbitrarily
high class and order, respectively. This illustrates geometrically
how, in the setting of knot concordance, the Vassiliev (isotopy)
invariants `collapse' to the Arf invariant.

\end{abstract}

\maketitle

\section{Introduction}

This paper gives a geometric characterization of the Arf invariant
of a knotted circle in the 3--sphere that is related to recent
developments in knot theory \cite{T1}. Conant and Teichner have
shown \cite{CT1} that the Vassiliev finite type filtration
\cite{B} on isotopy classes of knots corresponds to a geometric
equivalence relation called {\em 3--dimensional capped
grope-cobordism} and that this equivalence relation is generated
by certain simplified {\em half-gropes} in $S^3$. Gropes (see
\cite{T2}) are 2--complexes built by gluing together embedded
surfaces and in this setting the Vassiliev degree corresponds to a
measure of grope complexity called {\em class}, which counts the
layers of attached surfaces. The Arf invariant $\Arf(k)\in\Z_2$ of
a knot $k\subset S^3$ is the mod 2 reduction of the lowest degree
nontrivial finite type knot invariant (the degree two coefficient
of the Conway polynomial of $k$) and a result of Ng \cite{Ng}
says that $\Arf(k)$ is the only finite type invariant of $k$ up to
{\em concordance}, that is, up to bordism of $k$ by an embedded
annulus in the product $S^3\times I$ of the 3--sphere with an
interval. Gropes have been extensively studied in 4--dimensional
topology (e.g.~\cite{FQ,FT,Kr,KQ,KrT}) and are closely related to {\em Whitney towers} which
measure the failure of the Whitney move in terms of intersections
among higher order Whitney disks \cite{COT,S1,ST1,ST2}.

The following theorem illustrates geometrically how, in the
setting of knot concordance, the Vassiliev isotopy invariants
``collapse'' to the Arf invariant and the failure of the Whitney
move can be ``pushed out'' to arbitrarily high order Whitney
disks:
\begin{thm}\label{arf-thm}
For knots $k_0$ and $k_1$ in $S^3$ the following are equivalent:
\begin{enumerate}

\item $\mathrm{Arf}(k_0)=\mathrm{Arf}(k_1)$,

\item $k_0$ and $k_1$ cobound a properly embedded class $n$ half-grope in $S^3\times I$ for all $n\in\N$,

\item  $k_0$ and $k_1$ cobound a properly immersed annulus in $S^3\times I$ admitting
an order $n$ Whitney tower for all $n\in\N$.

\end{enumerate}

\end{thm}
It follows from Theorem~\ref{arf-thm} and a result in \cite{S1}
(which describes how to convert gropes into half-gropes) that a
knot in the 3-sphere has trivial Arf invariant if and only if it
bounds embedded gropes (not necessarily half-gropes) of
arbitrarily high class in the 4-ball.

Definitions of Whitney towers, $\Arf(k)$ and half-gropes will be
given in sections \ref{sec:w-tower}, \ref{sec:Arf} and
\ref{sec:half-grope-simple-tower}, respectively.
Theorem~\ref{arf-thm} will be proved constructively by exploiting
the flexibility of the Whitney towers in (iii). The infinite
cyclic Vassiliev (isotopy) invariant which lifts the Arf invariant
can be interpreted as the obstruction to ``pushing down'' this
construction into the 3--sphere.

\subsubsection*{Remark} The equivalence of (i) and (ii) also follows (somewhat indirectly)
from results in \cite{CT1}, \cite{CT2} and \cite{COT} (Proposition
3.8 in \cite{CT2}). \vspace{.07 in}

It should be mentioned that slight variations of the bordism
equivalence relations of {\em class $n$ grope concordance} and
{\em order $n$ Whitney concordance} suggested by (ii) and (iii)
are highly nontrivial. For instance, when using ``height'' instead
of class and order to measure complexity \cite{COT}, Cochran and
Teichner have used von Neumann $\rho$-invariants to show that the
associated filtration on grope (and Whitney) concordance classes
of knots is non-trivial for all $n$ \cite{CoT}. Also, when
working with links rather than knots, the leading term of the tree
part of the Kontsevich integral (equivalently \cite{HM}, Milnor's
$\mu$-invariants) gives obstructions to increasing the class
(resp. order) of a grope concordance (resp. Whitney concordance)
of the link \cite{ST2}.

\subsection*{Simple Whitney towers} The essential arguments in the proof of Theorem~\ref{arf-thm}
are contained in two lemmas. The first describes a close
relationship between half-gropes and certain {\em simple} Whitney
towers (\ref{simple-w-tower-subsec}) both of which are geometric
analogues of simple (right- or left-normed) commutators in a group
\cite{MKS}.

\begin{lem}\label{Grope-tower-lem}

Let $L$ be a link in the boundary of a simply connected
4--manifold. Then $L$ is the boundary of the bottom stage of a
properly embedded half-grope of class $n$ if and only if $L$ is
the boundary of a properly immersed planar surface which admits a
simple Whitney tower of order $n-1$.

\end{lem}

\subsubsection*{Remark} A much more general relation between class $n$ gropes
and order $n-1$ Whitney towers in (not
necessarily simply connected) 4--manifolds is described in
\cite{S1}. \vspace{.07 in}

Whitney towers are of interest in their own right, in part because
an order $n$ Whitney tower comes equipped with an $n$th order
geometric intersection obstruction to the existence of an
$(n+1)$th order tower which is related to Milnor's
($\overline{\mu}$) link invariants and the Kontsevich
integral \cite{ST2}, and is conjectured to generalize to give homotopy
invariants of immersed surfaces in arbitrary 4-manifolds
\cite{ST3}. (See also the more recent papers \cite{CST0,CST2,CST5,WTCCL}.) The next lemma illustrates how in the present
setting (of knots in a simply connected manifold) the obstruction
theory collapses after order 2.

\begin{lem}\label{w-tower-lem}

A properly immersed connected surface in a simply connected
4--manifold admitting an order 2 Whitney tower admits an order $n$
simple Whitney tower for all $n$.

\end{lem}

The connectivity conditions in Lemma~\ref{w-tower-lem} are
crucial. For instance (as explained in \cite{ST2}), in the setting
of {\em link} concordance, a first non-vanishing term of Vassiliev
degree $n$ in the tree part of the Kontsevich integral is an
obstruction to building a Whitney tower of order $n$ on a
collection of immersed disks in the 4--ball bounded by the link
components in $S^3$, and the above mentioned higher order Whitney
tower intersection obstruction is conjectured to be (highly)
non-trivial for connected surfaces in \emph{non}-simply connected
4-manifolds.

\subsection*{Outline} Whitney towers are defined in Section~\ref{sec:w-tower}, which also describes
the basic geometric manipulations of immersed surfaces in
4--manifolds that will be used throughout. The Arf invariant is
defined in Section~\ref{sec:Arf}. Lemma~\ref{Grope-tower-lem} is
proved in Section~\ref{sec:half-grope-simple-tower}, which also
contains definitions of half-gropes and simple Whitney towers. A
proof of Lemma~\ref{w-tower-lem} is given in
Section~\ref{sec:w-tower-lem-proof} and the proof of
Theorem~\ref{arf-thm} is assembled in
Section~\ref{sec:arf-thm-proof}. All manifolds are assumed smooth
and oriented.

%\tableofcontents

%%%%%%%%%%%%%%%%%%%%%%%%%%%%%%%%%%%%%%%%%%%%%%%%%%%%%%%%%%%%%%%%%%%%%%%%%%%%%%%%%%%%%%%%%%%%%

\section{Whitney towers}\label{sec:w-tower}
Whitney towers are introduced in this section, along with some
fundamental techniques from the theory of immersed surfaces in
4--manifolds. More information about Whitney towers can be found
in \cite{COT}, \cite{CST}, \cite{S1}, \cite{ST1} and \cite{ST2}.
For more details on surfaces in 4--manifolds the reader is
referred to \cite{FQ}.

\begin{figure}[ht!]
        \centerline{\includegraphics[scale=.6]{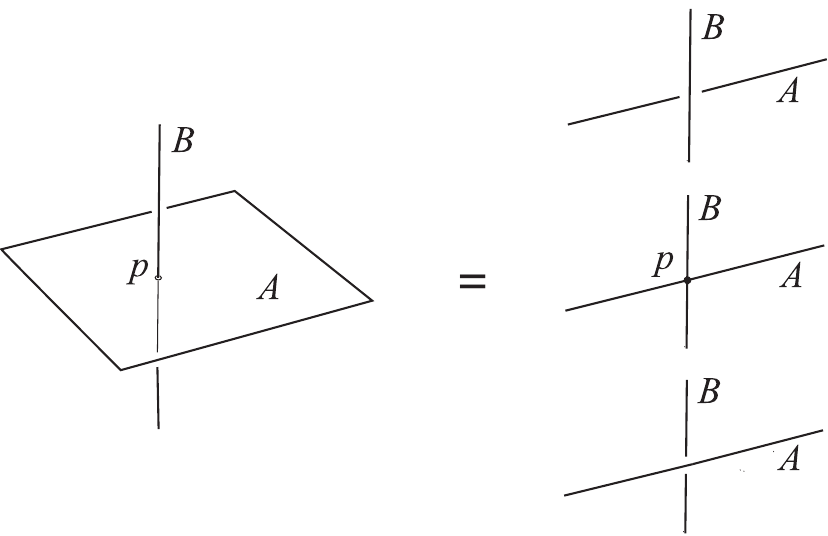}}
        \caption{A generic transverse intersection point $p$
         between surfaces $A$ and $B$ in a 4--manifold.}
        \label{int-point-fig}

\end{figure}

It will be convenient to illustrate surfaces locally in 4--space
by picturing 3--dimensional slices in which a surface may appear
either in the ``present slice'' or as an arc which extends into
neighboring slices; surfaces may also appear as a ``movie of
arcs'' in a sequence of 3--dimensional slices
(Figure~\ref{int-point-fig}).

\subsection{Whitney disks}\label{w-disk-subsec}
\begin{figure}[ht!]
        \centerline{\includegraphics[scale=.35]{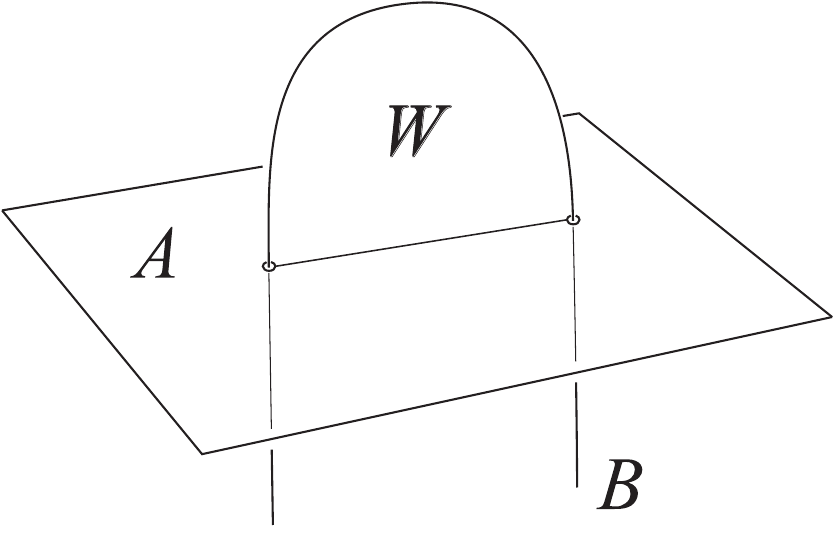}}
        \caption{An embedded Whitney disk $W$ pairing intersections between surfaces $A$ and $B$.}
         \label{W-disk-fig}

\end{figure}

In a simply connected 4--manifold $X$, two intersection points
between oriented connected surfaces $A$ and $B$ are a called a
{\em cancelling pair} if they have opposite signs (via the usual
sign convention that compares the orientations of the surfaces at
an intersection point with the orientation of the ambient
manifold). Such a cancelling pair $p$ and $q$ in $A\cap B$ can be
{\em paired by a Whitney disk} as follows: The union of an arc
$\alpha$ from $p$ to $q$ in $A$ and an arc $\beta$ from $q$ to $p$
in $B$ forms a loop in $X$ which bounds an immersed 2--disk $W$
meeting $A$ and $B$ along $\partial W$ in the standard way. Such a
$W$ is a {\em Whitney disk} pairing $p$ and $q$. (An embedded
Whitney disk is shown in Figure~\ref{W-disk-fig}.)
The normal disk bundle $\nu W$ of $W$ in $X$ pulls back to a
trivial $D^2$-bundle over the pre-image of $W$ (which is
contractible). The restriction of $\nu W$ to $\partial W$ has a
canonical 1--dimensional sub-bundle $\nu_{\partial}$ which
restricts along $\alpha$ to the normal bundle of $\alpha$ in $A$
and restricts along $\beta$ to the orthogonal complement (in $\nu
W$) of the normal bundle of $\beta$ in $B$. Since $p$ and $q$ have
opposite signs, $\nu_{\partial}$ is a trivial $I$-bundle over
$\partial W$. The obstruction to extending a non-vanishing section
of $\nu_{\partial}$ to a non-vanishing section of $\nu W$ is an
element of $\pi_1 SO(2)\cong \Z$. If this obstruction vanishes
then $W$ is said to be {\em framed}.

\subsection{Definition of a Whitney tower}\label{subsec:w-tower}
A Whitney disk can be used to eliminate its cancelling pair of
intersection points via a {\em Whitney move} (a motion of one of
the sheets guided by the Whitney disk), as introduced by Whitney
for higher dimensional manifolds immersed in Euclidean space
\cite{Wh}. In the present 4--dimensional setting, a Whitney move
will create new intersections if the interior of the Whitney disk
has any ``higher order'' intersections with sheets of surfaces or
Whitney disks (or if the Whitney disk is not framed), whereas in
dimensions greater than 4, such higher order intersection points
can eliminated by general position. The following notion of a {\em
Whitney tower} filters the condition that a properly immersed
surface in a 4--manifold is homotopic (rel boundary) to an
embedding.
\begin{figure}[ht!]
         \centerline{\includegraphics[scale=.65]{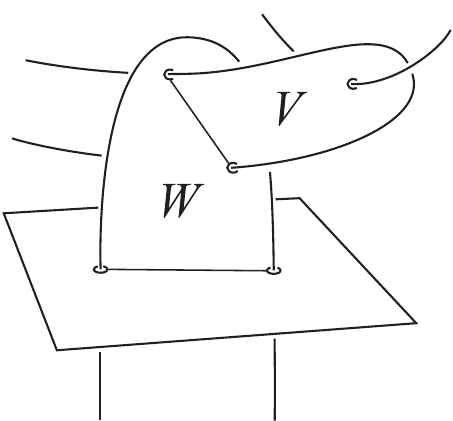}}
         \caption{Part of an order $n$ Whitney tower.
         The Whitney disk $V$ contains an unpaired intersection
         point which must be of order greater than or equal to $n$.}
         \label{W-tower-part-fig}

\end{figure}
\begin{defi}\label{w-tower-defi}\mbox{}

\begin{itemize}
\item A {\em surface of order 0} in a 4--manifold $X$
is a properly immersed surface (boundary embedded in the boundary
of $X$ and interior immersed in the interior of $X$). A {\em
Whitney tower of order 0} in $X$ is a collection of order 0
surfaces.
\item The {\em order of a (transverse) intersection point} between a surface of order $n$ and a
surface of order $m$ is $n+m$.
\item The {\em order of a Whitney disk} is $(n+1)$ if it pairs intersection points of order $n$.
\item For $n\geq 0$,
a {\em Whitney tower of order $(n+1)$}  is a Whitney tower $\cW$
of order $n$ together with (framed) Whitney disks pairing all
order $n$ intersection points of $\cW$. (These top order disks are
allowed to intersect each other as well as lower order surfaces.)
\end{itemize}

All Whitney disks in a Whitney tower are oriented (arbitrarily)
and are required to have disjointly embedded boundaries.

If $A$ is a properly immersed surface in a 4--manifold and there
exists an order $n$ Whitney tower containing $A$ as its order 0
surface, then $A$ is said to {\em admit} an order $n$ Whitney
tower.

\end{defi}
If a Whitney tower of order $n$ has no intersection points of
order greater than or equal to $n$, then the Whitney disks can be
used to guide a regular homotopy (rel boundary) of the order $0$
surfaces to an embedding.

\subsection{Modifying Whitney
disks}\label{modifying-w-disk-subsec} There are several moves that
allow for controlled modification of Whitney towers. Since the
moves are supported in a neighborhood of an arc or a point they
commute with each other and can be iterated disjointly arbitrarily
many times.
\begin{figure}[ht!]
         \centerline{\includegraphics[scale=.60]{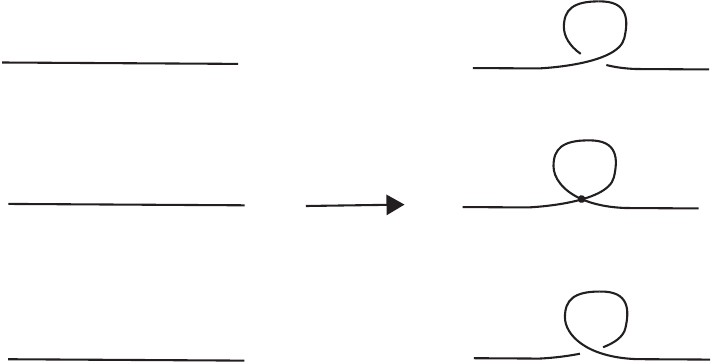}}
         \caption{A local cusp homotopy.}
         \label{cusp-fig}

\end{figure}

The first two {\em twisting moves} change the framing obstruction of a Whitney
disk.
\subsubsection{Interior twisting.}\label{interior-twist-subsec}
Introducing an self-intersection in the interior $\int W$ of $W$
by a {\em cusp homotopy} (see Figure~\ref{cusp-fig}) changes the
framing obstruction by $\pm2$ as can be seen by counting the intersections
between a local kink and its parallel push off.
\begin{figure}[ht!]
         \centerline{\includegraphics[scale=.5]{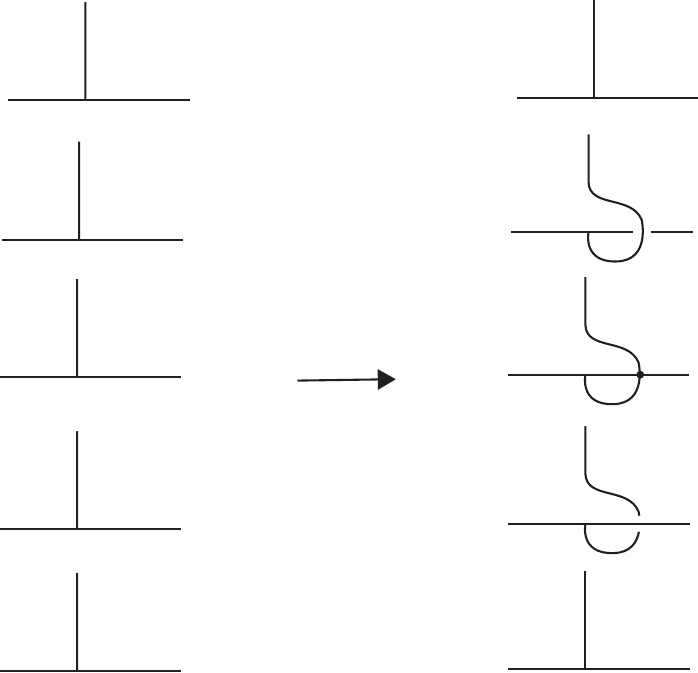}}
        \caption{Boundary twisting a Whitney disk.}
         \label{boundary-twist-fig}

\end{figure}
\subsubsection{Boundary twisting.}\label{boundary-twist-subsec}
Introducing a {\em boundary twist} by changing a collar of $W$
near a point in $\partial W$ (as in
Figure~\ref{boundary-twist-fig}) changes the framing obstruction by $\pm 1$
and creates an intersection between $\int W$ and the sheet
containing $\partial W$.

The next two moves do not affect framing obstructions but can be used to make
Whitney disks disjointly embedded.
\subsubsection{Boundary push-off.}\label{boundary-push}
\begin{figure}[ht!]
         \centerline{\includegraphics[scale=.4]{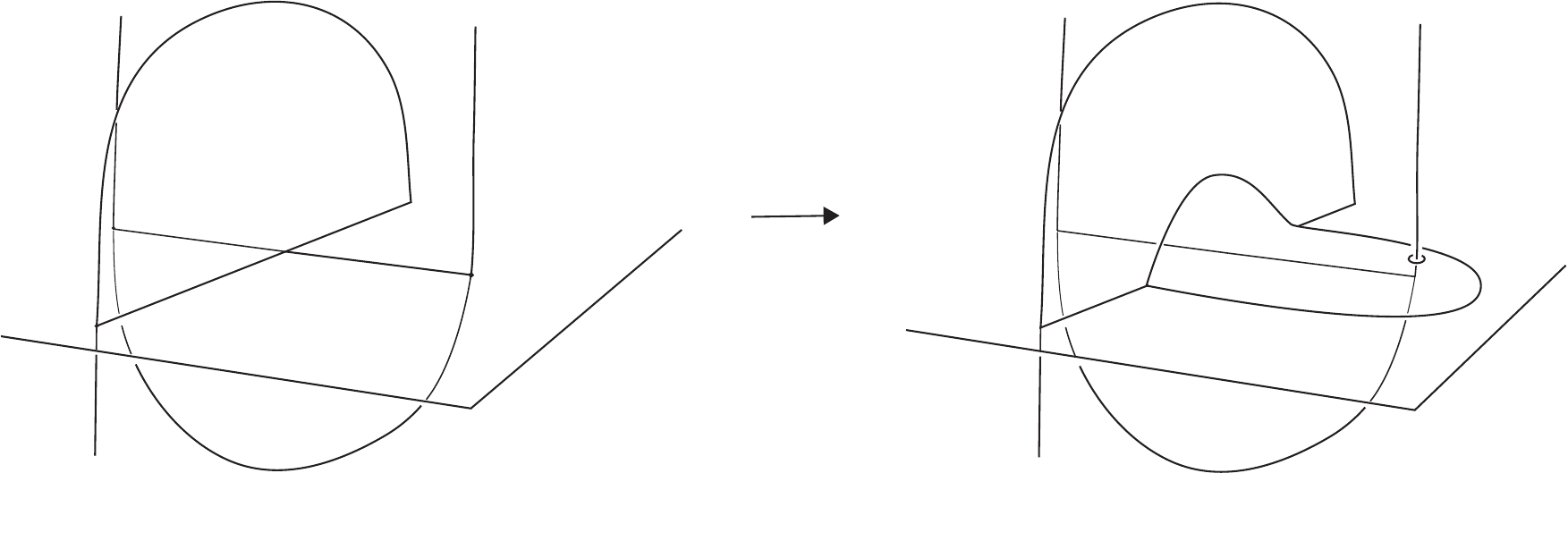}}
        \caption{Boundary push-off.}
         \label{boundary-push-fig}

\end{figure}

Intersections or self-intersections between boundaries of Whitney
disks can always be eliminated by a regular homotopy in a collar
at the cost of creating an intersection between the interior of a
Whitney disk and a surface sheet (Figure~\ref{boundary-push-fig}).

\subsubsection{Pushing down an intersection point.}\label{push-down}
\begin{figure}[ht!]
         \centerline{\includegraphics[scale=.4]{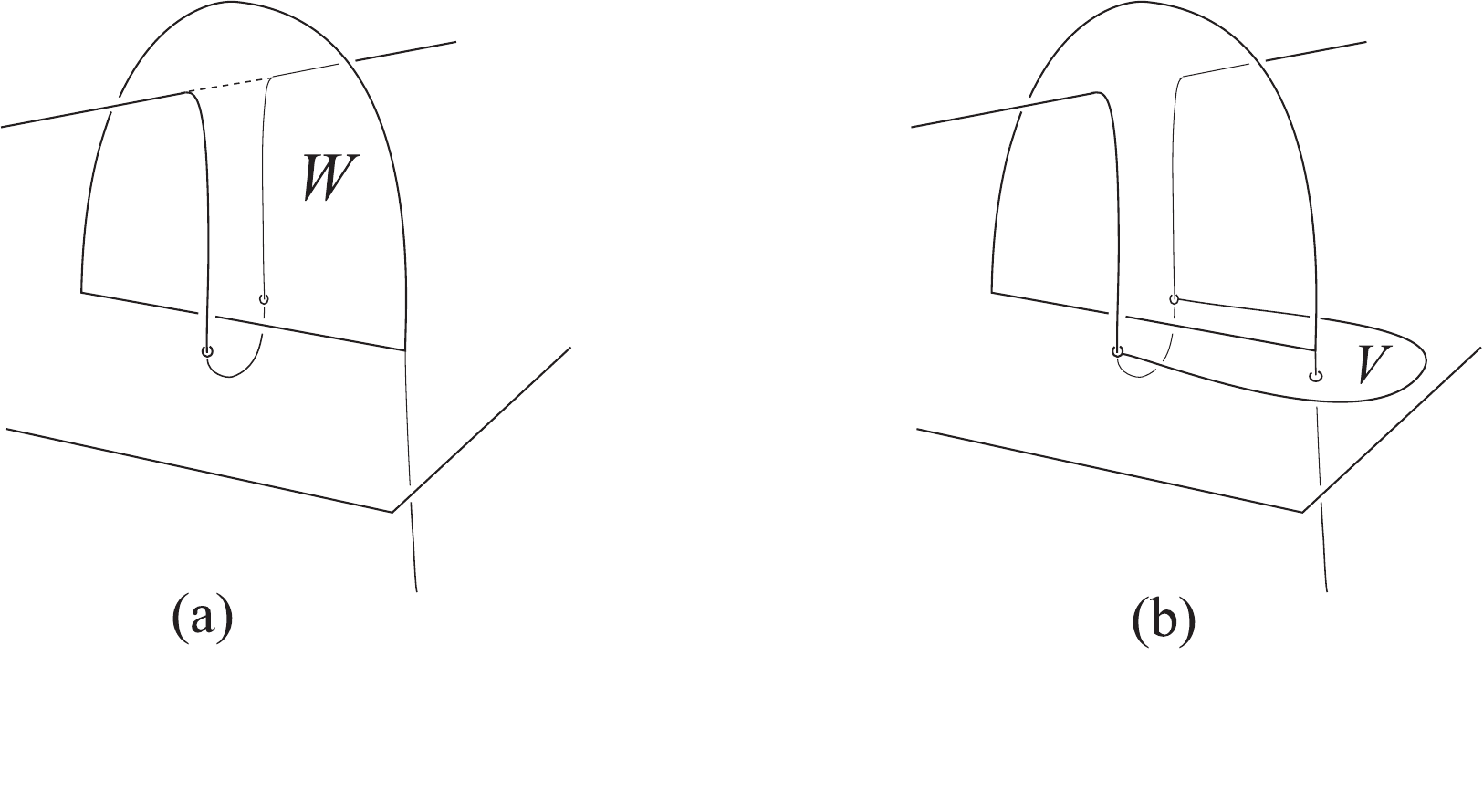}}
         \caption{(a) Pushing down an intersection between the interior of
         a Whitney disk $W$ and a sheet of a surface.
         (b) A Whitney disk $V$ for the cancelling pair created by pushing down
         (with $\partial V\cap\partial W=\emptyset$).}
         \label{push-down-fig}

\end{figure}

An intersection point between $\int W$ and any surface $S$ can be
removed by a {\em finger move} on $S$, a homotopy of $S$ supported
in a neighborhood of an arc, in this case an arc in $W$ from the
intersection to a point in $\partial W$
(Figure~\ref{push-down-fig}(a)). Such a finger move is called
``pushing $S$ down'' into either sheet and creates two new
cancelling pairs of intersection points between $S$ and the sheet.
This move can also be used to remove interior self-intersections
of $W$.
Note that the newly created cancelling pair can be paired by an
embedded Whitney disk $V$ whose boundary is disjoint from $W$ by
applying the boundary push-off move to the obvious small embedded
Whitney disk near the cancelling pair. The interior of $V$ has a
single intersection point with the sheet that was not ``pushed
into'' as illustrated in Figure~\ref{push-down-fig}(b).

\subsection{Order 1 towers for knots}
Applying the moves of \ref{modifying-w-disk-subsec} yields the
following lemma, which will be used in the definition of the Arf
invariant in the next section.
\begin{lem}\label{1-tower-lem}
Any knot $k$ in $S^3$ bounds a properly immersed 2-disk in
$S^3\times I$ admitting a Whitney tower of order 1.
\end{lem}
\begin{proof}
A finite number of crossing changes, leading from $k$ to the
unknot, describes a properly immersed 2-disk $D$ in $S^3\times I$
(with the unknot capped off by an embedded disk). Fixing
orientations, the signs of the self-intersections of $D$
correspond to the signs of the crossing changes and, after
introducing trivial crossing changes (if necessary), the
self-intersections of $D$ can be made to occur in cancelling pairs
of order~0 intersections which are paired by order~1 Whitney disks
as in Section~\ref{w-disk-subsec}. By applying boundary twists
(\ref{boundary-twist-subsec}) and boundary push-off
(\ref{boundary-push}) as needed, it can be arranged that the
Whitney disks are framed with disjointly embedded boundaries.
\end{proof}

\section{The Arf invariant}\label{sec:Arf}
In \cite{Ro}, Robertello used Kervaire and Milnor's generalization
\cite{KM} of Rochlin's Theorem to define a $\Z_2$-valued
concordance invariant of a knot in $S^3$, and showed that it was
equal to the Arf invariant of a quadratic enhancement of (the mod
2 reduction of) the Seifert form. (The Arf invariant of a
non-degenerate quadratic form is defined to be 0 (resp. 1) if a
majority of elements are taken to 0 (resp. 1).) This knot
invariant has numerous characterizations, all of which are
commonly referred to as the Arf invariant (of a knot). The
following geometric definition, which we have translated into the
language of Whitney towers, is due to Matsumoto \cite{Ma} using
Freedman and Kirby's geometric proof of Rochlin's Theorem
\cite{FK}.

\begin{defi}\label{arf-defi}
For a knot $k$ in $S^3$, let $D$ be any properly immersed 2-disk
immersed in $S^3\times I$ admitting a Whitney tower $\cW$ of
order~1, with $k=\partial D\subset S^3\times \{0\}$. Define the
{\em Arf invariant} of $k$, $\Arf(k)\in\Z_2$, to be the number
(modulo 2) of order~1 intersection points in $\cW$.
\end{defi}
By Lemma~\ref{1-tower-lem}, such a $\cW$ always exists. A direct
combinatorial proof that this definition of $\Arf(k)$ is
well-defined can be found in $10.8$ of \cite{FQ}.
\begin{figure}[ht!]
        \centerline{\includegraphics[scale=.45]{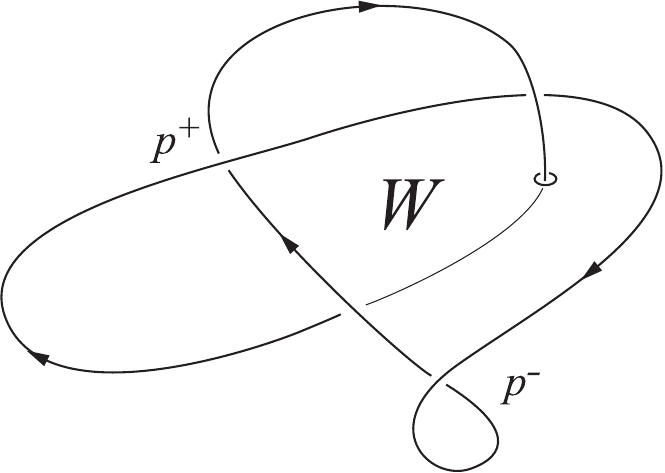}}
         \caption{The trefoil has nontrivial Arf invariant.}
         \label{trefoil-arf-fig}

\end{figure}
\subsection{The trefoil.}
That the trefoil knot has nontrivial Arf invariant can be seen in
Figure~\ref{trefoil-arf-fig}, where changing the crossings
labelled $p^+$ and $p^-$ creates a cancelling pair of order 0
intersections in a null-homotopy $D$. This cancelling pair has a
framed embedded Whitney disk $W$ which intersects $D$ in a single
order~1 intersection point.

\section{Half-gropes and simple Whitney towers}\label{sec:half-grope-simple-tower}
In this section a proof of Lemma~\ref{Grope-tower-lem} is given
after first defining half-gropes and simple Whitney towers, two
geometric analogues of a simple (right- or left-normed) commutator
of elements in a group \cite{MKS}. It should be noted that the
fact that we are working with {\em half}-gropes and {\em simple}
Whitney towers is crucial in the below proof of
Lemma~\ref{Grope-tower-lem}. In the setting of general gropes and
Whitney towers, showing the correspondence between class and order
involves more subtle geometric constructions (see \cite{S1}).
Basic operations on gropes (surgery, etc.) used in this section
are described in detail in \cite{FQ}.
\begin{figure}[ht!]
         \centerline{\includegraphics[scale=.4]{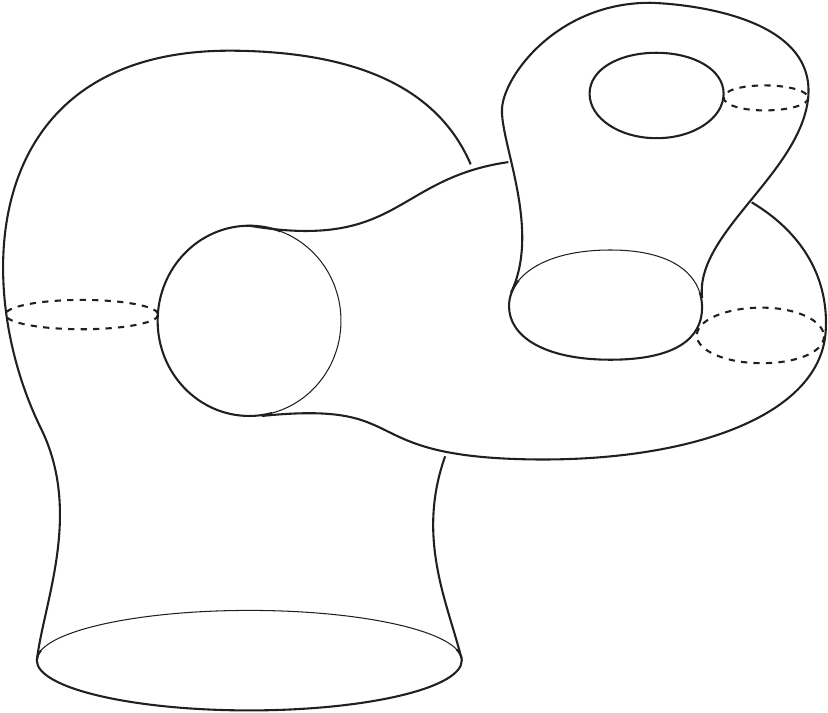}}
         \caption{A half-grope of class 4.}
         \label{half-grope-fig}

\end{figure}

\subsection{Half-gropes}\label{half-grope-subsec}
In general, gropes are 2-complexes consisting of surfaces joined
along certain essential curves (see e.g. \cite{T2,CT1,FQ}). Requiring that the curves form a ``half-basis'' yields
the ``half-gropes'':
\begin{defi}\label{grope-defi}

A {\em half-grope of class $2$} is a compact connected orientable
surface $A$ with a single boundary circle. To form a {\em
half-grope of class $n>2$}, start with an orientable surface $A$
with a single boundary circle and choose a symplectic basis
$\{a_i,b_i\}$, that is, the $a_i$ and $b_i$ are embedded curves
which represent a basis for $H_1(A)$ and the only intersections
among them occur when $a_i$ meets $b_i$ in a single point. Now
attach half-gropes of class $n-1$ along their boundary circles to
a $1/2$-symplectic basis for $A$, i.e., a maximal pairwise
disjoint subset of $\{a_i,b_i\}$, for instance $\{a_i\}$.

The surface $A$ is called the {\em bottom stage} of the half-grope
and the boundary circle of $A$ is the {\em boundary} of the
half-grope. Half-gropes with more than one boundary component are
formed by removing disks from the bottom stage of a half-grope.
The attached punctured surfaces are also referred to as {\em
(higher) stages}. The basis curves that do not have higher stages
attached to them are the {\em tips} of the half-grope.

A half-grope $H$ is {\em properly embedded} in a 4--manifold $X$
if the boundary of $H$ is embedded in $\partial X$ and the rest of
$H$ is embedded in $\int X$. It is also required that $H$ satisfy
the following normal {\em framing condition}: A regular
neighborhood of $H$ in $X$ must factor as a standard embedding of
$H$ into 3--space followed by taking the product with an interval.

\end{defi}

\subsection{Simple Whitney towers}\label{simple-w-tower-subsec}
A Whitney tower is {\em simple} if all of its Whitney disks have
disjointly embedded interiors. Thus, every intersection point of
order $m$ in a simple Whitney tower is an intersection between a
surface of order $0$ and a surface of order $m$.

\subsection{Proof of Lemma~\ref{Grope-tower-lem}}\label{Grope-tower-lem-proof}
Let $H$ be a half-grope of class $n$ properly embedded in $X$ and
bounded by $L$. Since $X$ is simply connected, the tips of $H$
bound immersed 2-disks called {\em caps} and the plan is to create
the desired Whitney tower by surgering the caps. Each cap has a
normal framing obstruction determined by pushing its boundary along $H$ and
after boundary twisting the cap (just as in
\ref{boundary-twist-subsec}) this obstruction can be made to vanish.
We may also arrange, by repeatedly pushing down intersections as
for Whitney disks (\ref{push-down}), that the caps of $H$ are
disjointly embedded (except for the single boundary point
intersections between dual caps in the top stages) with interiors
disjoint from all stages of $H$ except for perhaps the bottom
stage.
\begin{figure}[ht!]
         \centerline{\includegraphics[scale=.45]{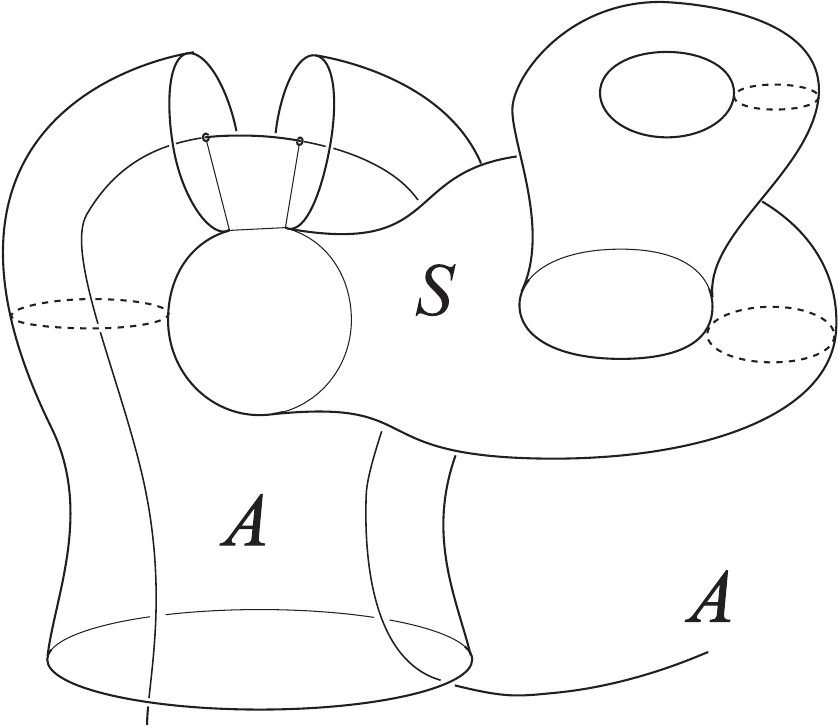}}
         \caption{Intersections created by surgering caps can be paired by Whitney disks
         constructed from surgered higher stages of the half-grope. Surgering a cap bounded by
         the dotted circle (tip) on $S$ creates a Whitney disk.}
         \label{surgered-grope-fig}
\end{figure}
Assume first that the bottom stage intersects the interior of each
cap in at most a single point. Let $A$ be the result of surgering
those caps attached to the bottom stage surface of $H$. Then $A$
is a properly immersed planar surface with self-intersection
points coming in cancelling pairs which were created by surgering
a cap whose interior intersected the bottom stage surface. Such a
cancelling pair has an embedded (first order) Whitney disk $W$
which is the union of a small band and the result of surgering the
caps on the next stage surface $S$, which was attached along the
dual curve to the boundary of the cap
(Figure~\ref{surgered-grope-fig}). The framing condition (in
Definition~\ref{grope-defi}) on the normal bundle of $H$ in $X$
ensures that $W$ is framed. The only possible
intersections between the interior of $W$ and anything else are
intersections with $A$ coming from intersections between $A$ and
the surgered caps on $S$, hence occur in cancelling pairs with an
embedded second order Whitney disk gotten by similarly surgering
the next surface stage. This construction terminates at the top
$(n-1)$th stage surfaces, where only a half-basis of caps are
surgered to make the order $(n-2)$ Whitney disks and the dual caps
(together with bands) form the order $(n-1)$ Whitney disks
yielding the desired simple Whitney tower.

The above assumption that each cap has at most a single interior
intersection with the bottom surface can always be arranged by
Krushkal's grope splitting technique \cite{Kr}; alternatively,
the above construction can still be carried out for ``near-by''
cancelling pairs created by surgering caps containing multiple
interior intersections by using parallel disjoint copies
(guaranteed by the normal framing condition) of the higher
surfaces stages of $H$ to build the higher order Whitney disks.
This completes one direction of the proof of
Lemma~\ref{Grope-tower-lem}.
\begin{figure}[ht!]
         \centerline{\includegraphics[scale=.5]{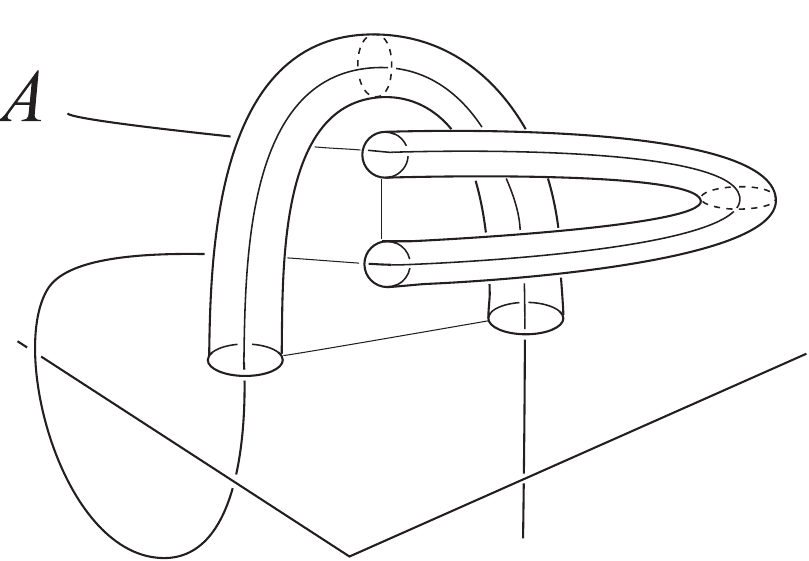}}
         \caption{A simple Whitney tower of order
         $n-1$ yields a half-grope of class $n$ by `tubing along $A$'.}
         \label{w-grope-fig}

\end{figure}
For the other direction, let $\cW$ be a simple order $n-1$ Whitney
tower on a properly immersed planar surface $A$ in $X$ bounded by
$L$. Since $\cW$ is simple, each of its Whitney disks has (at
least one) boundary arc lying in $A$ (the order 1 Whitney disks
have both boundary arcs on $A$). The desired half-grope $H$ of
class $n$ is constructed by ``tubing the Whitney disks of $\cW$
along $A$'' as illustrated in Figure~\ref{w-grope-fig}: More
specifically, let $W$ be any order $m$ Whitney disk ($1\leq m\leq
n-1$) for a pair of cancelling order $m-1$ intersections between
$A$ and an order $m-1$ Whitney disk $V$ (if $m=1$, then $V$ is
just the order zero surface $A$). Denote by $\partial_AW$ the part
of the boundary of $W$ that lies in $A$ (if $m=1$, choose a
boundary arc of $\partial W$). Using the boundary annulus of the
normal disk bundle to $A$ in $X$ restricted to $\partial_AW$ to
perform 0-surgery on $V$ eliminates the cancelling pair of
intersections between $A$ and $V$. If $m=n-1$, then $W$ is
discarded; if $1\leq m<n-1$, then $W$ (minus a small collar near
$\partial_AW$) becomes an ($m+1$)th stage surface of $H$ by
0-surgering the interior of $W$ to eliminate any intersections
with $A$.  Applying this construction to all the Whitney disks of
$\cW$ yields $H$, with the bottom stage surface of $H$ consisting
of 0-surgery (one for each first order Whitney disk) on $A$ and
each $m$th stage surface 0-surgery on an $(m-1)$th order Whitney
disk. The normal framing condition on $H$ is satisfied since all
the Whitney disks of $\cW$ were framed. \hfill$\square$

\section{Proof of Lemma~\ref{w-tower-lem}}\label{sec:w-tower-lem-proof}
The idea of the proof of Lemma~\ref{w-tower-lem} is that, for an
order $n$ simple Whitney tower whose order $0$ surface is {\em
connected}, any order $n$ intersection point can be cancelled by
boundary twisting its order $n$ Whitney disk into the order 0
surface (in a simply connected 4--manifold). The framing on the
Whitney disk can then be recovered by boundary twisting along the
other boundary arc of the Whitney disk, which only creates an
intersection point of order $n+(n-1)$, which is greater than $n$
for $n\geq 2$.

\begin{proof}
Let $A$ be a properly immersed connected surface in a simply
connected 4--manifold $X$ admitting an order 2 Whitney tower
$\cW$. A {\em simple} order 2 Whitney tower can be constructed
from $\cW$ by pushing down any intersections among its Whitney
disks into $A$: Pushing down (order 2) intersections among the
first order Whitney disks creates cancelling pairs of (order 1)
intersections which can be equipped with disjointly embedded
second order Whitney disks having a single interior (order 2)
intersection with $A$ (Figure~\ref{push-down-fig}). Pushing down
intersections among second order disks and between second and
first order disks creates cancelling pairs of second order
intersections between the second order disks and $A$.

Now assume inductively that $A$ admits a simple Whitney tower
$\cW_n$ of order $n\geq 2$. Since the interiors of all Whitney
disks in $\cW_n$ are disjointly embedded, the only possible
unpaired intersection points are $n$th order intersections between
$n$th order Whitney disks and $A$.

Let $p$ be such an intersection point between an $n$th order
Whitney disk $W$ and $A$. Since $\cW_n$ is simple, $W$ pairs
intersections between $A$ and an $(n-1)$th order Whitney disk
(recall $n\geq2$). By performing a boundary twist around the arc
of $\partial W$ that lies on $A$, we can create an intersection
point $q\in \int W\cap A $ of opposite sign as $p$. Since $A$ is
connected (and the 4--manifold is simply-connected), $p$ and $q$
can be paired by an $(n+1)$th order Whitney disk.
\begin{figure}[ht!]
         \centerline{\includegraphics[scale=.50]{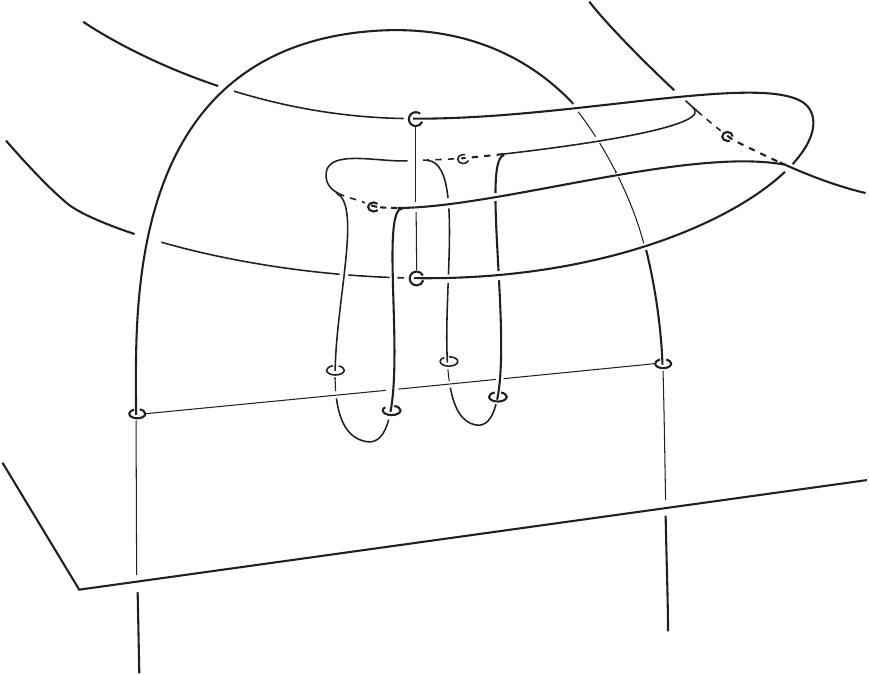}}
         \caption{}
         \label{multi-pushdown-fig}

\end{figure}
To restore the framing of $W$, perform a boundary twist along the
(other) arc of $\partial W$ that lies on the $(n-1)$th order
Whitney disk, creating a $(2n-1)$th order intersection point
between $W$ and the $(n-1)$th order Whitney disk. This $(2n-1)$th
order intersection point can be eliminated by repeatedly pushing
$\int W$ down into Whitney disks of lower order as in
Figure~\ref{multi-pushdown-fig} until eventually reaching $A$,
where $2^{(n-2)}$ cancelling pairs of order $n$ intersections
between $\int W$ and $A$ will be created. These cancelling pairs
admit disjointly embedded order $n$ Whitney disks (parallel copies
of the Whitney disk $V$ pictured in
Figure~\ref{push-down-fig}(b)), each having a single order~$n+1$
intersection with $A$.

Since this modification of $\cW_n$ takes place in a neighborhood
of a 1-complex, it may be repeated (in disjoint neighborhoods)
until all order $n$ intersections are paired by order $(n+1)$
Whitney disks. The boundaries of these $(n+1)$th order Whitney
disks can be made disjointly embedded (and disjoint from all other
Whitney disk boundaries) by applying boundary push-off moves
(\ref{boundary-push}). Finally, intersections between any Whitney
disk and the $(n+1)$th order Whitney disks can be eliminated by
repeatedly pushing the $(n+1)$th order Whitney disks down (as in
Figure~\ref{multi-pushdown-fig}) until they only intersect $A$,
yielding a simple order $(n+1)$ Whitney tower.
\end{proof}

%%%%%%%%%%%%%%%%%%%%%%%%%%%%%%%%%%%%%%%%%%%%%%%%%%%%%%%%%%%%%%%%%%%%%%%%%%%%%%%%%%%%%%%%%%%%%

%--------------------------------------------------------------------------------------------

\section{Proof of Theorem~\ref{arf-thm}}\label{sec:arf-thm-proof}
\subsection*{{\boldmath $\mathrm{(iii)}\Rightarrow\mathrm{(i)}$}}
Let $A$ be an annulus admitting an order $n\geq 2$ Whitney tower
$\cW$ as in (iii). Then any 2-disk $D_0$ admitting an order 1
Whitney tower $\cW_0$ in $S^3\times I$ and bounded by $k_0$ can be
extended by $A$ to a 2-disk $D_1=A\cup D_0$ in $S^3\times I$
admitting an order 1 Whitney tower $\cW_1=\cW \cup \cW_0$ and
bounded by $k_1$. Since $n\geq 2$, all order~1 intersection points
in $\cW$ occur in cancelling pairs, so $\cW_1$ has the same number
(modulo 2) of order 1 intersection points as $\cW_0$ and
$\Arf(k_1)=\Arf(k_0)$.
\begin{figure}[ht!]
         \centerline{\includegraphics[scale=.55]{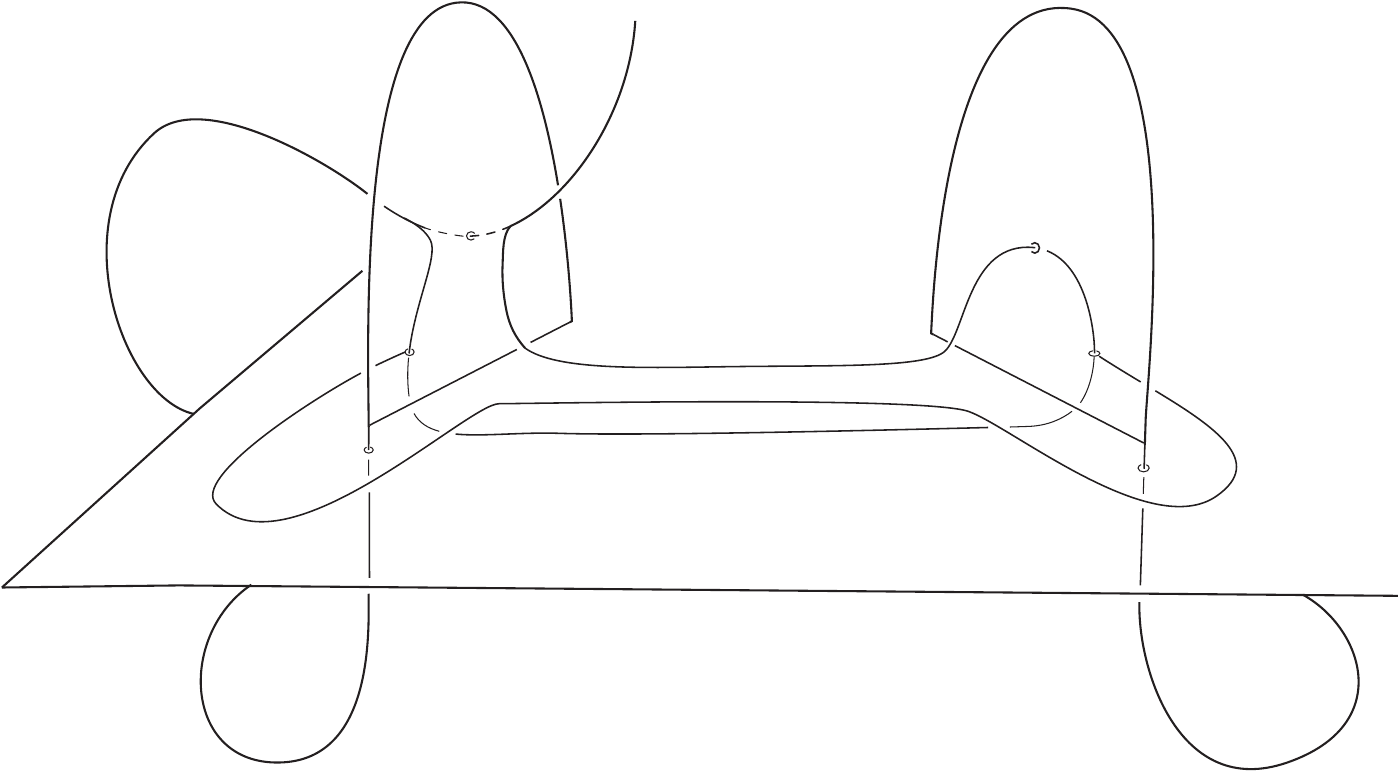}}
         \caption{}
         \label{Arf-move-fig}

\end{figure}
\subsection*{{\boldmath $\mathrm{(i)}\Rightarrow\mathrm{(iii)}$}}
Let $D_i$, $i=0,1$, be disjoint properly immersed 2-disks in
$S^3\times I$, where $I=[0,1]$, bounded by the knots $k_i \subset
S^3\times\{i\}$ and admitting order 1 Whitney towers $\cW_i$. We
may assume that the $\cW_i$ are disjoint, so the $D_i$ can be
tubed together by a thin embedded annulus to get an annulus $A$
co-bounded by the $k_i$ such that $A$ admits an order 1 Whitney
tower $\cW$ whose Whitney disks are just the union of the Whitney
disks in $\cW_i$. The assumption that $\Arf(k_0)=\Arf(k_1)$ means
that $\cW$ has an even number of order~1 intersection points. By
using the move illustrated in Figure~\ref{Arf-move-fig} (details
in \cite{Y}, also \cite{ST1}), which does not affect framing obstructions, it
can be arranged that each Whitney disk contains an even number of
order~1 intersection points. After introducing an even number of
(like-signed) boundary twists (\ref{boundary-twist-subsec}) on
each Whitney disk, the order~1 intersection points on each Whitney
disk occur in cancelling pairs (have opposite signs) admitting
second order Whitney disks. If $2m$ boundary twists were done on a
first order Whitney disk $W$, then the framing of $W$ can be
recovered by performing $m$ interior twists
(\ref{interior-twist-subsec}), which create only second order
intersection points (self-intersections of $W$). Having thus far
constructed a second order Whitney tower, the proof of the
implication $\mathrm{(i)}\Rightarrow\mathrm{(iii)}$ is completed
by Lemma~\ref{w-tower-lem}.

\hfill
\subsection*{{\boldmath $\mathrm{(ii)}\Leftrightarrow \mathrm{(iii)}$}}
By Lemma~\ref{w-tower-lem}, we may assume that the Whitney tower
in (iii) is simple, thus the proof of
$\mathrm{(ii)}\Leftrightarrow \mathrm{(iii)}$ follows from
Lemma~\ref{Grope-tower-lem}.\hfill$\square$

%=============================================================================

\end{document}